\documentclass[11pt]{amsart}
\usepackage{mathrsfs}
\usepackage{amssymb,latexsym}
\usepackage{pstricks,color}

 \numberwithin{equation}{section}

\newtheorem{theorem}{Theorem}[section]

\newtheorem{corollary}[theorem]{Corollary}

\newtheorem{lemma}[theorem]{Lemma}

\def\qed{\hfill $\Box$}
\def\pf{\noindent {\it Proof.} }

\title{Enumerating several statistics of $r$-Colored Dyck paths with no $\mathbf{dd}$-steps having the same colors }

\begin{document}
\maketitle
\begin{center}
Yidong Sun\footnote{Corresponding author: Yidong Sun}, Jinyi Wang$^2$ and Xinyu Wang$^3$

School of Science, Dalian Maritime University, 116026 Dalian, P.R. China\\[5pt]

{\it Emails: $^1$sydmath@dlmu.edu.cn, $^2$wangjinyi@dlmu.edu.cn, $^3$wxymath@dlmu.edu.cn}

\end{center}\vskip0.2cm

\subsection*{Abstract}  An $r$-colored Dyck path is a Dyck path with all $\mathbf{d}$-steps having one of $r$ colors in $[r]=\{1, 2, \dots, r\}$. In this paper, we consider several statistics on the set $\mathcal{A}_{n,0}^{(r)}$ of $r$-colored Dyck paths of length $2n$ with no two consecutive $\mathbf{d}$-steps having the same colors. Precisely, the paper studies the statistics ``number of points" at level $\ell$, ``number of $\mathbf{u}$-steps" at level $\ell+1$, ``number of peaks" at level $\ell+1$ and ``number of $\mathbf{udu}$-steps" on the set $\mathcal{A}_{n,0}^{(r)}$. The counting formulas of the first three statistics are established by Riordan arrays related to $S(a,b; x)$, the weighted generating function of $(a,b)$-Schr\"{o}der paths. By a useful and surprising relations satisfied by $S(a,b; x)$, several identities related to these counting formulas are also described.

\medskip

{\bf Keywords}: Dyck path; Riordan array; $(a, b)$-Schr\"{o}der path; $(a, b)$-Schr\"{o}der number.

\noindent {\sc 2000 Mathematics Subject Classification}: Primary
05A15; Secondary 05A10, 05A19.

{\bf \section{ Introduction } }

A {\it Dyck path} of length $2n$ is a lattice path from $(0, 0)$ to $(2n, 0)$ in the first quadrant of the $XY$-plane and consists of $n$ up-steps $\mathbf{u}=(1, 1)$ and $n$ down-steps $\mathbf{d}=(1, -1)$. See \cite[p.204]{StanleyEC} and \cite{Deutsch99}. Let $\mathcal{D}_n$ be the set of Dyck paths of length $2n$ and $\mathcal{D}=\bigcup_{n\geq 0}\mathcal{D}_n$. It is well known \cite{Stanley} that $|\mathcal{D}_n|=C_n=\frac{1}{n+1}\binom{2n}{n}$, the $n$th Catalan number. An {\it $r$-colored Dyck path} is a Dyck path with all $\mathbf{d}$-steps having one of $r$ colors in $[r]=\{1, 2, \dots, r\}$.

A {\it Motzkin path} of length $n$ is a lattice path from $(0, 0)$ to $(n, 0)$ in the first quadrant of the $XY$-plane and consists of up-steps $\mathbf{u}=(1, 1)$, down-steps $\mathbf{d}=(1, -1)$ and horizontal-steps $\mathbf{h}=(1, 0)$.
A {\it Schr\"{o}der path} of length $2n$ is a path from $(0, 0)$ to $(2n, 0)$ in the first quadrant of the XY-plane that consists of up steps $\mathbf{u}=(1, 1)$, down steps $\mathbf{d}=(1, -1)$ and horizontal steps $\mathbf{H}=(2, 0)$.

An {\it $(a,b)$-Dyck path} is a weighted Dyck path with $\mathbf{u}$-steps weighted by $1$, $\mathbf{d}$-steps in peaks (called $\mathbf{ud}$-steps) weighted by $a$ and other $\mathbf{d}$-steps weighted by $b$. An {\it $(a,b)$-Motzkin path} of length $n$ is a weighted Motzkin path with $\mathbf{u}$-steps weighted by $1$, $\mathbf{h}$-steps weighted by $a$ and $\mathbf{d}$-steps weighted by $b$. An {\it $(a,b)$-Schr\"{o}der path} is a weighted Schr\"{o}der path such that the $\mathbf{u}$-steps, $\mathbf{H}$-steps, and $\mathbf{d}$-steps are weighted respectively by $1, a$ and $b$.
The {\it weight} of a path $\mathbf{P}$, denoted by $w(\mathbf{P})$, is the product of the weight of each step of $\mathbf{P}$. The weight of a set $\mathcal{A}$ of any weighted paths, denoted by $w(\mathcal{A})$, is the sum of the total weights of all paths in $\mathcal{A}$.

Let $\mathcal{C}_n(a,b), \mathcal{M}_n(a,b)$ and $\mathcal{S}_n(a,b)$ be respectively the sets of $(a,b)$-Dyck paths of length $2n$, $(a,b)$-Motzkin paths of length $n$ and $(a,b)$-Schr\"{o}der paths of length $2n$. Let $C_n(a,b)$, $M_n(a,b)$ and $S_n(a,b)$ be their weights with $C_0(a,b)=M_0(a,b)=S_0(a,b)=1$ respectively.
It is known \cite{ChenPan} that $C_n(a,b), M_n(a,b)$ and $S_n(a,b)$ have the explicit formulas:
\begin{eqnarray}
C_n(a,b)  \hskip-.22cm &=&\hskip-.22cm  \sum_{k=1}^{n}\frac{1}{n}\binom{n}{k-1}\binom{n}{k}a^{k}b^{n-k},   \label{eqn 1.1} \\
M_n(a,b)  \hskip-.22cm &=&\hskip-.22cm  \sum_{k=0}^{n}\binom{n}{2k}C_ka^{n-2k}b^{k},        \nonumber          \\
S_n(a,b)  \hskip-.22cm &=&\hskip-.22cm  \sum_{k=0}^{n}\binom{n+k}{2k}C_ka^{n-k}b^{k},  \label{eqn 1.2}
\end{eqnarray}
and their generating functions
\begin{eqnarray}
C(a,b; x)  \hskip-.22cm &=&\hskip-.22cm \sum_{n=0}^{\infty}C_n(a,b)x^n=\frac{1-(a-b)x-\sqrt{(1-(a-b)x)^2-4bx}}{2bx}, \nonumber \\
M(a,b; x)  \hskip-.22cm &=&\hskip-.22cm \sum_{n=0}^{\infty}M_n(a,b)x^n=\frac{1-ax-\sqrt{(1-ax)^2-4bx^2}}{2bx^2},     \nonumber \\
S(a,b; x)  \hskip-.22cm &=&\hskip-.22cm \sum_{n=0}^{\infty}S_n(a,b)x^n=\frac{1-ax-\sqrt{(1-ax)^2-4bx}}{2bx},  \label{eqn 1.3}
\end{eqnarray}
which satisfy the functional equations:
\begin{eqnarray}
C(a,b; x)  \hskip-.22cm &=&\hskip-.22cm 1+(a-b)xC(a,b; x)+bxC(a,b; x)^2, \nonumber\\
M(a,b; x)  \hskip-.22cm &=&\hskip-.22cm 1+axM(a,b; x)+bx^2M(a,b; x)^2,   \nonumber\\
S(a,b; x)  \hskip-.22cm &=&\hskip-.22cm 1+axS(a,b; x)+bxS(a,b; x)^2. \label{eqn 1.4}
\end{eqnarray}

There are closely relations between $C_n(a,b)$, $M_n(a,b)$ and $S_n(a,b)$. Exactly, Chen and Pan \cite{ChenPan} derived the following equivalent relations
\begin{eqnarray}\label{eqn 1.5}
S_n(a,b)=C_n(a+b,b)=(a+b)M_{n-1}(a+2b,(a+b)b)
\end{eqnarray}
for $n\geq 1$ and provided some combinatorial proofs by weighted Dyck, Motzkin and Schr\"{o}der paths. Other combinatorial proofs of (\ref{eqn 1.5}) have been given by \cite{SunWang} using weighted generalized Motzkin paths.

When $a=b=1$, $C_n(a,b), M_n(a,b)$ and $S_n(a,b)$ reduce to the Catalan number $C_n$ \cite{Deutsch99, Stanley}, the Motzkin number $M_n$ \cite{Aigner, BarPinSpr, DonShap} and the large Schr\"{o}der number $S_n$ \cite{BoninShap, Sloane} respectively. Write $C(x), M(x), S(x)$ instead of $C(a,b; x), M(a,b; x)$, $S(a,b; x)$ respectively when $a=b=1$.

Let $\varepsilon$ be the empty path, i.e., the path of length $0$. If both $\mathbf{P}_1$ and $\mathbf{P}_2$ are of same type (Dyck, Motzkin or Schr\"{o}der paths), then we define $\mathbf{P}_1\mathbf{P}_2$ as the concatenation of $\mathbf{P}_1$ and $\mathbf{P}_2$. Any non-empty Dyck path $\mathbf{P}$ has a unique first return decomposition \cite{BarKirPet1,Deutsch99} of the form $\mathbf{P}= \mathbf{u}\alpha\mathbf{d}\beta$, where $\alpha$ and $\beta$ are
two Dyck paths in $\mathcal{D}$.

For any path such as Dyck, Motzkin or Schr\"{o}der path, a point of a path with ordinate $\ell$ is said to be at {\it level} $\ell$. A step of a path is said to be at level $\ell$ if the ordinate of its endpoint is $\ell$. By a {\it return step} we mean a $\mathbf{d}$-step at level $0$. A {\it peak (valley)} in a path is an occurrence of $\mathbf{ud}$ ($\mathbf{du}$). By the {\it level of a peak (valley)} we mean the level of the intersection point of its two steps. Dyck paths that have exactly one return step are said to be {\it primitive}.

In the literature, many papers deal with the enumeration of Dyck paths or their restricted classes according to different statistics, e.g. {\it peak}, {\it valley}, {\it pyramid}, {\it water capacity}, {\it number of $\mathbf{udu}'s$} \cite{Asakly, BleBre, Czabarka1, Czabarka2, DenSim, Deutsch99, Elizalde, FerMun, FloJunRam, FloRam19, FloRam20, ManSap, Man2002, Man2006, MerSprVer, PanSap, PeaWoan, SapTasTsi, Sun, SunJia}, and very important statistics such as {\it area, dinv, bounce} of Dyck paths, see for example \cite{DAddIraWyn, GarXin, Haglund, XinZhang19, XinZhang23} and references therein. Other papers deal with Motzkin paths, Schr\"{o}der paths or various generalized paths using similar methods \cite{BarPinSpr, BarKirPet2, BreMav, WPCheng, DonShap, ManSchSun, SapTsi, SunWang, SunZhao, WagProd}. For instance, the Shapiro's Catalan triangle \cite{Shap76} $\mathcal{B}=(B_{n,k})_{n\geq k\geq 0}$ with $B_{n,k}=\frac{k+1}{n+1}\binom{2n+2}{n-k}$ counts the number of points at level $k$ of Dyck paths of length $2n$ \cite{WPCheng}, or counts the number of peaks at level $k+1$ of Dyck paths of length $2n+2$ \cite{WPCheng}. Note that this triangle forms a well-kown Riordan array \cite{SunMaCat, SunMaMot} $(C(x)^2, xC(x)^2)$. Another example is that the general entry $A_{n,k}=\frac{2k+3}{2n+3}\binom{2n+3}{n-k}$ of the Riordan array \cite{SunMaCat, SunMaMot} $(C(x)^3, xC(x)^2)$ counts the number of $\mathbf{u}$-steps at level $k+1$ of Dyck paths of length $2n+2$ \cite{WPCheng}.

Recall that {\it Riordan array} \cite{Rogers78, ShapB, Shap18, ShapGet, Sprug} is an infinite lower triangular matrix $\mathscr{D}=(d_{n,k})_{n,k \in \mathbb{N}}$ such that its $k$-th column has generating function $d(x)h(x)^k$, where $d(x)$ and $h(x)$ are formal power series with $d(0)\neq 0$ and $h(0)=0$. That is, the general term of $\mathscr{D}$ is $d_{n,k}=[x^n]d(x)h(x)^k$, where $[x^n]$ is the coefficient operator. The matrix $\mathscr{D}$ corresponding to the pair $d(x)$ and $h(x)$ is denoted by $(d(x),h(x))$. A Riordan array $\mathscr{D}=(d(x),h(x))$ is {\it proper}, if $h'(0) \neq 0$ additionally. Suppose we multiply the matrix $\mathscr{D} = (d(x), h(x))$ by a column vector $(a_0, a_1, \dots)^T$ and get a column vector $(b_0, b_1, \dots)^T$. Let $A(x)$ and $B(x)$ be the generating functions for the sequences $(a_0, a_1, \dots)^T$ and $(b_0, b_1, \dots)^T$ respectively. Then it follows easily that
$$B(x)=(d(x), h(x))A(x)=d(x)A(h(x)).$$

Riordan arrays play an important role in combinatorics in terms of obtaining combinatorial identities and have applications in many fields of mathematics. See, for example, \cite{Barry17, Barry21, ChenLiangWang, He18, MeRogSprugVer, ShapGet, ShapSpruBarChHe, Slowik, Sprug, SunSun, SunMaMot, YangYang, ZhangZhao} and references therein.

In this paper, we consider several statistics on $r$-colored Dyck paths with no two consecutive $\mathbf{d}$-steps (called $\mathbf{dd}$-steps) having the same colors. Precisely, the next section dedicates on $r$-colored Dyck paths with no $\mathbf{dd}$-steps having the same colors and introduces a useful lemma related to $S(a,b; x)$, the weighted generating function of $(a,b)$-Schr\"{o}der paths. The third section studies the statistic ``number of points" at level $\ell$. The fourth section concentrates on the statistic ``number of $\mathbf{u}$-steps" at level $\ell+1$. The fifth section devotes to the statistic ``number of peaks ($\mathbf{ud}$-steps)" at level $\ell+1$. The last section investigates the statistic ``number of $\mathbf{udu}$-steps".  The counting formulas of the first three statistics are established by Riordan arrays related to $S(a,b; x)$, several identities related to these counting formulas are also described.

\vskip0.5cm

\section{ $r$-Colored Dyck paths }

An $r$-colored Dyck path is a Dyck path with all $\mathbf{d}$-steps having one of $r$ colors in $[r]=\{1, 2, \dots, r\}$.
Let $\mathcal{A}_{n,k}^{(r)}$ denote the set of $r$-colored Dyck paths $\mathbf{P}$ of length $2n$ such that $\mathbf{P}$ exactly has $k$ $\mathbf{dd}$-steps with the same color. Set $\mathcal{A}_{k}^{(r)}=\bigcup_{n\geq k}\mathcal{A}_{n,k}^{(r)}$. Let $A_{n,k}^{(r)}$ be the cardinality of $\mathcal{A}_{n,k}^{(r)}$ and $A_r(x, y)$ be the bivariate generating function of $A_{n,k}^{(r)}$, i.e.,
$$A_r(x, y)=\sum_{n=0}^{\infty}A_{n}^{(r)}(y)x^{n}=\sum_{n=0}^{\infty}\sum_{k=0}^{n}A_{n,k}^{(r)}y^{k}x^{n}. $$

Denote by $A_{r,i}(x, y)$ to be the bivariate generating function of the number of paths $\mathbf{P}\in \mathcal{A}_{n,k}^{(r)}$ according to $n$ and $k$ such that the last return step of $\mathbf{P}$ is colored by $i$. Such nonempty paths $\mathbf{P}$ can be written as $\mathbf{P}=\mathbf{P}_1\mathbf{u}\mathbf{P}_2\mathbf{d}_i$, according to the last return decomposition, where $\mathbf{P}_1\in \mathcal{A}_{n_1, k_1}^{(r)}, \mathbf{P}_2\in \mathcal{A}_{n_2, k_2}^{(r)}$ for some nonnegative integers $n_1, n_2, k_1$ and $k_2$ such that $n_1+n_2=n-1$ and $k_1+k_2=k-1\ (k_1+k_2=k)$ if the color of the last return step of $\mathbf{P}_2$ is (not) $i$ , and $\mathbf{d}_i$ denotes the last return step of $\mathbf{P}$ with the color $i$. Moreover, once the last return step of $\mathbf{P}_2$ is colored by $i$, it contributes one additional $\mathbf{dd}$-step with the same color to $\mathbf{P}$. Then we deduce the following relations:
\begin{eqnarray*}
A_r(x, y)     \hskip-.25cm &=& \hskip-.25cm 1+ \sum_{i=1}^{r}A_{r,i}(x, y),  \\
A_{r,i}(x, y) \hskip-.25cm &=& \hskip-.25cm x\big(A_r(x, y)+(y-1)A_{r,i}(x, y)\big)A_r(x, y).
\end{eqnarray*}

These produce a recursive relation for $A_r(x, y)$, namely,
\begin{eqnarray*}
A_r(x, y) \hskip-.25cm &=& \hskip-.25cm 1+ x(1-y)A_r(x, y) + x(y+r-1)A_r(x, y) ^2.
\end{eqnarray*}
By solving it, we have
\begin{eqnarray*}
A_r(x, y) \hskip-.25cm &=& \hskip-.25cm  \frac{1-(1-y)x-\sqrt{(1-(1-y)x)^2-4(y+r-1)x}}{2(y+r-1)x},
\end{eqnarray*}
which generates that $A_{0}^{(r)}(y)=1$ and
\begin{eqnarray*}
A_{n}^{(r)}(y) \hskip-.25cm &=& \hskip-.25cm \sum_{k=0}^{n}A_{n,k}^{(r)}y^{k}=C_n(r, y+r-1)=S_n(1-y, y+r-1)    \\
\hskip-.25cm &=& \hskip-.25cm \sum_{\ell=1}^{n}\frac{1}{n}\binom{n}{\ell-1}\binom{n}{\ell}r^{\ell}(y+r-1)^{n-\ell}
\end{eqnarray*}
for $n\geq 1$. Taking the coefficient of $y^k$, we obtain the following result.

\begin{theorem}\label{theo 2.1}
The set of $r$-colored Dyck paths of length $2n$ with exact $k$ $\mathbf{dd}$-steps having the same color is counted by
\begin{eqnarray}\label{eqn 2.1}
A_{n,k}^{(r)} \hskip-.25cm &=& \hskip-.25cm \sum_{\ell=1}^{n-k}\frac{1}{n}\binom{n}{\ell-1}\binom{n}{\ell}\binom{n-\ell}{k}r^{\ell}(r-1)^{n-\ell-k}.
\end{eqnarray}
In the $k=0$ case $A_{n, 0}^{(r)}=S_n(1, r-1)$ counts the number of $r$-colored Dyck paths of length $2n$ with no $\mathbf{dd}$-steps having the same color. More precisely, the number of $2$-colored Dyck paths of length $2n$ with no $\mathbf{dd}$-steps having the same color is enumerated by the large Schr\"{o}der number $S_n$ \cite{Sloane}.
\end{theorem}

The formula (\ref{eqn 2.1}) has a direct combinatorial interpretation as follows. It is well known that the Narayana number $\frac{1}{n}\binom{n}{\ell-1}\binom{n}{\ell}$ counts the Dyck paths $\mathbf{P}$ of length $2n$ with exact $\ell$ peaks. If one colors the $\ell$ $\mathbf{d}$-steps in peaks by one of the $r$ colors, there are $r^{\ell}$ ways. Then one has $\binom{n-\ell}{k}$ ways to choose $k$ $\mathbf{d}$-steps from the left $n-\ell$ $\mathbf{d}$-steps such that each of the $k$ $\mathbf{d}$-steps has the same color as the $\mathbf{d}$-step right ahead of it. Finally, each of the rest $n-\ell-k$ $\mathbf{d}$-steps has different colors as the $\mathbf{d}$-step immediately ahead of it, there are $(r-1)^{n-\ell-k}$ ways to do this. In conclusion, we derive (\ref{eqn 2.1}).

Next, we present the Lagrange inversion formula \cite{Gessel} as a necessary tool and provide a useful lemma related to $S(a,b; x)$ as follows.

{\bf The Lagrange inversion formula}: Let $u(t)$ be a formal power series not involving $x$ and $u(0)\neq 0$. Then there is a unique formal power series
$f(x)$ such that $f(x)=xu(f(x))$, and for any formal power series $\Psi(t)$ not involving $x$ and any integer $n$ we have
\begin{eqnarray*}
[x^{n}]\Psi(f(x)) \hskip-.22cm &=&\hskip-.22cm \frac{1}{n}[t^{n-1}]\Psi'(t)u(t)^{n}
\end{eqnarray*}
for all $n\geq 1$, where $\Psi'(t)$ is the derivative of $\Psi(t)$ with respect to $t$.

\begin{lemma}\label{lemma 2.2}
Let $Z_1(a, b; x)=S(a,b; -x)S(a,b; \frac{b}{a}xS(a,b; -x)^2)$ for $a\neq 0$. Then $Z_1(a, b; x)$ satisfies the relation
\begin{eqnarray}\label{eqn 2.2}
Z_1(a, b; x) \hskip-.25cm &=& \hskip-.25cm 1-axZ_1(a, b; x)+\frac{b^2}{a}xZ_1(a, b; x)^2.
\end{eqnarray}
Moreover, let
$$Z_{m+1}(a, b; x)=Z_m(a, b; x)Z_m(a,b; \frac{b^{2^{m}}}{a^{2^{m}}}xZ_m(a,b; x)^2)$$
with $Z_0(a, b; x)=S(-a, b; x)$ for $a\neq 0$ and $m\geq 0$. Then $Z_m(a, b; x)$ satisfies the relation
\begin{eqnarray}\label{eqn 2.3}
Z_m(a, b; x) \hskip-.25cm &=& \hskip-.25cm 1-axZ_m(a, b; x)+\frac{b^{2^m}}{a^{2^m-1}}xZ_m(a, b; x)^2
\end{eqnarray}
and has the explicit formula
\begin{eqnarray}\label{eqn 2.4}
Z_m(a, b; x) \hskip-.25cm &=& \hskip-.25cm S(-a, \frac{b^{2^m}}{a^{2^m-1}}; x) = \frac{1+ax-\sqrt{(1+ax)^2-4\frac{b^{2^m}}{a^{2^m-1}}x}}{2\frac{b^{2^m}}{a^{2^m-1}}x}.
\end{eqnarray}
Specially, $Z_m(a, b; x)=Z_m(a, -b; x) $ for $m\geq 1$ and $Z_m(a, a; x)=S(-a, a; x)=1$ for $m\geq 0$.
\end{lemma}

\pf Setting $T:=S(a,b; -x)$, by (\ref{eqn 1.4}), we have
\begin{eqnarray*}
Z_1(a, b; x) \hskip-.22cm &=&\hskip-.22cm TS(a,b; \frac{b}{a}xT^2) \\
 \hskip-.22cm &=&\hskip-.22cm  T\Big(1+bxT^2S(a,b; \frac{b}{a}xT^2)+\frac{b^2}{a}xT^2S(a,b; \frac{b}{a}xT^2)^2\Big) \\
 \hskip-.22cm &=&\hskip-.22cm  T+ bxT^2Z_1(a, b; x) + \frac{b^2}{a}xTZ_1(a, b; x)^2   \\
 \hskip-.22cm &=&\hskip-.22cm  T+ (1-axT-T)Z_1(a, b; x) + \frac{b^2}{a}xTZ_1(a, b; x)^2,
\end{eqnarray*}
which produces that
\begin{eqnarray*}
T\Big(1-ax Z_1(a, b; x)-Z_1(a, b; x) + \frac{b^2}{a}xZ_1(a, b; x)^2\Big) \hskip-.22cm &=&\hskip-.22cm 0.
\end{eqnarray*}
Since $T\neq 0$, we immediately obtain (\ref{eqn 2.2}).

If setting $Z_0(a, b; x)=S(-a, b; x)$, we verify similarly that
$$Z_{1}(a, b; x)=Z_0(a, b; x)Z_0(a,b; \frac{b}{a}xZ_0(a,b; x)^2) $$
also satisfies the relation (\ref{eqn 2.2}).

Similar to the proof of (\ref{eqn 2.2}), by induction on $m$, we get (\ref{eqn 2.3}). The detail is left to the interested readers. Specially, by (\ref{eqn 1.4}), it is easy to verify that (\ref{eqn 2.4}) holds, and (\ref{eqn 2.4}) implies that $Z_m(a, b; x)=Z_m(a, -b; x) $ for $m\geq 1$ and $Z_m(a, a; x)=S(-a, a; x)=1$. \qed

\begin{corollary}\label{coro 2.3}
For any integers $n, m\geq 0$ and any parameters $a, b$ with $a\neq 0$, we obtain
\begin{eqnarray*}
S_n\big(-a, \frac{b^{2^{m+1}}}{a^{2^{m+1}-1}}\big) \hskip-.22cm &=&\hskip-.22cm \sum_{\ell=0}^{n}\sum_{j=0}^{n-\ell}(-1)^{j}\frac{2\ell+1}{2n-j+1}\binom{2n-j+1}{n-\ell} \binom{n-\ell}{j}  \\
\hskip-.22cm & &\hskip-.22cm \hskip2cm \cdot  S_{\ell}\big(-a, \frac{b^{2^{m}}}{a^{2^{m}-1}}\big)a^{n-\ell-2^{m}(n-j)}b^{2^{m}(n-j)}.
\end{eqnarray*}
\end{corollary}

\pf Setting $Y_m=Z_{m}(a, b; x)-1$, by (\ref{eqn 2.3}), we deduce
\begin{eqnarray*}
Y_m=x(1+Y_m)\Big(\frac{b^{2^m}}{a^{2^m-1}}(1+Y_m)-a\Big).
\end{eqnarray*}

Taking the coefficient of $x^{n}$ in $Z_{m+1}(a, b; x)$ and using the Lagrange inversion formula to $Y_m$, one has
\begin{eqnarray*}
\lefteqn{S_n\big(-a, \frac{b^{2^{m+1}}}{a^{2^{m+1}-1}}\big) = [x^n]Z_{m+1}(a, b; x) } \\
\hskip-.22cm &=&\hskip-.22cm [x^n]Z_m(a, b; x)Z_m(a,b; \frac{b^{2^{m}}}{a^{2^{m}}}xZ_m(a,b; x)^2)  \\
\hskip-.22cm &=&\hskip-.22cm [x^n]\sum_{\ell=0}^{\infty}S_{\ell}\big(-a, \frac{b^{2^{m}}}{a^{2^{m}-1}}\big)\frac{b^{\ell2^{m}}}{a^{\ell2^{m}}}x^{\ell}Z_m(a, b; x)^{2\ell+1}     \\
\hskip-.22cm &=&\hskip-.22cm  \sum_{\ell=0}^{n}S_{\ell}\big(-a, \frac{b^{2^{m}}}{a^{2^{m}-1}}\big)\frac{b^{\ell2^{m}}}{a^{\ell2^{m}}}[x^{n-\ell}]Z_m(a, b; x)^{2\ell+1}               \\
\hskip-.22cm &=&\hskip-.22cm  \sum_{\ell=0}^{n}S_{\ell}\big(-a, \frac{b^{2^{m}}}{a^{2^{m}-1}}\big)\frac{b^{\ell2^{m}}}{a^{\ell2^{m}}}\frac{2\ell+1}{n-\ell}
[Y_m^{n-\ell-1}](1+Y_m)^{n+\ell}\Big(\frac{b^{2^m}}{a^{2^m-1}}(1+Y_m)-a\Big)^{n-\ell}            \\
\hskip-.22cm &=&\hskip-.22cm \sum_{\ell=0}^{n}\sum_{j=0}^{n-\ell}(-1)^{j}\frac{2\ell+1}{2n-j+1}\binom{2n-j+1}{n-\ell} \binom{n-\ell}{j}  \\
\hskip-.22cm & &\hskip-.22cm \hskip2cm\cdot  S_{\ell}\big(-a, \frac{b^{2^{m}}}{a^{2^{m}-1}}\big)a^{n-\ell-2^{m}(n-j)}b^{2^{m}(n-j)},
\end{eqnarray*}
which produces the result. \qed

In the sequel, we have main interest in the set $\mathcal{A}_{n,0}^{(r)}$ of $r$-colored Dyck paths of length $2n$ with no $\mathbf{dd}$-steps having the same colors. For convenience, we set $S_r(x)=\sum_{n\geq 0} S_n^{(r)}x^n$, where
$$ S_r(x):=A_r(x, 0)=S(1, r-1; x)=\frac{1-x-\sqrt{(1-x)^2-4(r-1)x}}{2(r-1)x}. $$
Then $S_n^{(r)}=A_{n,0}^{(r)}=S_n(1, r-1)$ is the number of $r$-colored Dyck paths of length $2n$ with no $\mathbf{dd}$-steps having the same colors. The first values of $S_n^{(r)}$ for $2\leq r\leq 5$ and $0\leq n\leq 8$ are illustrated in Table 2.1.
\begin{center}
\begin{eqnarray*}
\begin{array}{c|ccccccccc|c}\hline
  n        & 0      & 1      & 2       & 3       & 4      & 5      &   6       &   7        &   8          &   OEIS           \\\hline
S_n^{(2)}  & 1      & 2      & 6       & 22      & 90     & 394    & 1806      &  8558      &  41586       &   A006318 \  \mbox{\cite{Sloane}  }   \\\hline
S_n^{(3)}  & 1      & 3      & 15      & 93      & 645    & 4791   & 37275     &  299865    &  2474025     &   A103210 \  \mbox{\cite{Sloane}  }     \\\hline
S_n^{(4)}  & 1      & 4      & 28      & 244     & 2380   & 24868  & 272188    &  3080596   &  35758828    &   A103211 \  \mbox{\cite{Sloane}  }     \\\hline
S_n^{(5)}  & 1      & 5      & 45      & 505     & 6345   & 85405  & 1204245   &  17558705  &  262577745   &   A133305 \   \mbox{\cite{Sloane}  }     \\\hline
\end{array}
\end{eqnarray*}
Table 2.1. The first values of $S_{n}^{(r)}$.
\end{center}

Clearly, $S_r(x)$ satisfies the relation
$$S_r(x)=1+xS_r(x)+x(r-1)S_r(x)^2, $$
which implies that $S_2(x)=S(x)$, where $S(x)$ is the generating function of the large Schr\"{o}der numbers $S_n$ \cite{Sloane}. By Lemma 2.2 and (\ref{eqn 1.3}), we have
\begin{eqnarray}\label{eqn 2.5}
S_r(x)S_r(-x(r-1)S_r(x)^2)  \hskip-.25cm &=& \hskip-.25cm S(-1, (r-1)^2; -x)= S(1, -(r-1)^2; x).
\end{eqnarray}
By (\ref{eqn 1.1})-(\ref{eqn 1.3}) and (\ref{eqn 2.5}), similar to the proof of Corollary 2.3, one gets the identity,
\begin{eqnarray*}
\lefteqn{\sum_{k=0}^{n}(-1)^{k}\binom{n+k}{2k}C_k(r-1)^{2k}}\\
\hskip-.22cm &=&\hskip-.22cm \sum_{\ell=0}^{n}\sum_{j=0}^{n-\ell}(-1)^{\ell}\frac{2\ell+1}{2n-j+1}\binom{2n-j+1}{n-\ell} \binom{n-\ell}{j}S_{\ell}^{(r)}(r-1)^{n-j}.
\end{eqnarray*}
The $r=2$ case in (\ref{eqn 2.5}) indicates that
\begin{eqnarray*}
S(x)S(-xS(x)^2)  \hskip-.25cm &=& \hskip-.25cm 1
\end{eqnarray*}
which is an interesting relation related to the generating function $S(x)$ of the large Schr\"{o}der numbers $S_n$.

\section{The statistic ``number of points" at level $\ell$ on $\mathcal{A}_{n,0}^{(r)}$ }
Let $P_{n, \ell}^{(r)}$ denote the total number of points at level $\ell$ in all $r$-colored Dyck paths in $\mathcal{A}_{n,0}^{(r)}$ ($\mathcal{A}_{n,0}^{(r)}$ is the set of all $r$-colored Dyck paths of length $2n$ with no $\mathbf{dd}$-steps having the same color). The first values of $P_{n, \ell}^{(r)}$ for $r=2, 3$ are illustrated in Table 3.1 and Table 3.2 respectively.

\begin{center}
\begin{eqnarray*}
\begin{array}{c|ccccccc}\hline
n/\ell & 0      & 1      & 2       & 3       & 4    & 5    & 6    \\\hline
  0 & 1      &        &         &         &      &      &     \\
  1 & 4      & 2      &         &         &      &      &     \\
  2 & 16     & 12     & 2       &         &      &      &      \\
  3 & 68     & 64     & 20      & 2       &      &      &      \\
  4 & 304    & 332    & 144     & 28      &  2   &      &      \\
  5 & 1412   & 1712   & 916     & 256     & 36   &  2   &      \\
  6 & 6752   & 8844   & 5488    & 1948    & 400  &  44  & 2     \\\hline
\end{array}
\end{eqnarray*}
Table 3.1. The first values of $P_{n, \ell}^{(2)}$.
\end{center}

\begin{center}
\begin{eqnarray*}
\begin{array}{c|ccccccc}\hline
n/\ell & 0      & 1      & 2       & 3       & 4    & 5    & 6    \\\hline
  0 & 1      &        &         &         &       &      &        \\
  1 & 6      & 3      &         &         &       &      &        \\
  2 & 39     & 30     & 6       &         &       &      &        \\
  3 & 276    & 267    & 96      & 12      &       &      &        \\
  4 & 2073   & 2316   & 1128    & 264     & 24    &      &         \\
  5 & 16242  & 20031  & 11832   & 3876    & 672   & 48   &         \\
  6 & 131295 & 174018 & 117678  & 48000   & 11856 & 1632 & 96    \\\hline
\end{array}
\end{eqnarray*}
Table 3.2. The first values of $P_{n, \ell}^{(3)}$.
\end{center}

\begin{theorem}\label{theom 3.1.1}
For any integers $n\geq \ell\geq 0$ and $r\geq 2$, we deduce
\begin{eqnarray*}
P_{n, \ell}^{(r)}=\sum_{j=0}^{n-\ell}\binom{n-\ell}{j}\left\{\frac{2\ell}{(2n-j)(r-1)}\binom{2n-j}{n-\ell}+\frac{2\ell+2}{2n-j+2}\binom{2n-j+2}{n-\ell} \right\}(r-1)^{n-j}
\end{eqnarray*}
with $P_{0, 0}^{(r)}=1$. Moreover, $P_{n, \ell}^{(r)}+\frac{1}{r-1}\delta_{n, 0}$ is the $(n,\ell)$-entry of the Riordan array
$$\Big(S_r(x)^2+\frac{1}{r-1},\ x(r-1)S_r(x)^2\Big), $$
where $\delta_{n, 0}=1$ for $n=0$ and $0$ otherwise.
\end{theorem}

\pf If $\mathbf{D}$ is an empty path, then there is just one point at level $0$, namely, $P_{0, 0}^{(r)}=1$. For any $\mathbf{D}\in \mathcal{A}_{n+1,0}^{(r)}$, $\mathbf{D}$ can be partitioned into $\mathbf{D}=\mathbf{D}_1\mathbf{u}\mathbf{D}_2\mathbf{d}_j$ with $j\in [r]$,
$\mathbf{D}_1\in \mathcal{A}_{k,0}^{(r)}$ and $\mathbf{D}_2\in \mathcal{A}_{n-k,0}^{(r)}$ for some $0\leq k\leq n$.

For the number $P_{n+1, 0}^{(r)}$ of points at level $0$ in all $\mathbf{D}\in \mathcal{A}_{n+1,0}^{(r)}$, note that there are just $S_{n+1}^{(r)}$ number of $\mathbf{D}\in \mathcal{A}_{n+1,0}^{(r)}$ and each $\mathbf{D}$ only contributes one endpoint at $(2n+2, 0)$. In the decomposition $\mathbf{D}=\mathbf{D}_1\mathbf{u}\mathbf{D}_2\mathbf{d}_j$, all $\mathbf{D}_1\in \mathcal{A}_{k,0}^{(r)}$ have total $P_{k, 0}^{(r)}$ points at level $0$, there are $S_{n-k}^{(r)}$ choices for $\mathbf{D}_2\in \mathcal{A}_{n-k,0}^{(r)}$ and $r-1$ choices of colors for $\mathbf{d}_j$ if
$\mathbf{D}_2$ is nonempty and $r$ choices of colors for $\mathbf{d}_j$ if
$\mathbf{D}_2$ is empty. Then we built the following recurrence
\begin{eqnarray}\label{eqn 3.1}
P_{n+1, 0}^{(r)}  \hskip-.22cm &=&\hskip-.22cm (r-1)\sum_{k=0}^{n} P_{k, 0}^{(r)}S_{n-k}^{(r)}+S_{n+1}^{(r)}+P_{n, 0}^{(r)}.
\end{eqnarray}

For the number $P_{n+1, \ell}^{(r)}$ of points in all $\mathbf{D}\in \mathcal{A}_{n+1,0}^{(r)}$ at level $\ell\geq 1$, in the decomposition $\mathbf{D}=\mathbf{D}_1\mathbf{u}\mathbf{D}_2\mathbf{d}_j$, there are two cases to be considered.

Firstly, all $\mathbf{D}_1\in \mathcal{A}_{k,0}^{(r)}$ have total $P_{k, \ell}^{(r)}$ points at level $\ell$ for $\ell\leq k\leq n$, all $\mathbf{D}_2\in \mathcal{A}_{n-k,0}^{(r)}$ have $S_{n-k}^{(r)}$ choices and there are $r-1$ choices of colors for $\mathbf{d}_j$ if
$\mathbf{D}_2$ is nonempty and $r$ choices of colors for $\mathbf{d}_j$ if $\mathbf{D}_2$ is empty.

Secondly, all $\mathbf{D}_2\in \mathcal{A}_{n-k,0}^{(r)}$ have total $P_{n-k, \ell-1}^{(r)}$ points at level $\ell-1$ for $0\leq k\leq n-\ell+1$, all $\mathbf{D}_1\in \mathcal{A}_{k,0}^{(r)}$ have $S_{k}^{(r)}$ choices, when $\ell=1$ there are $r-1$ choices of colors for $\mathbf{d}_j$ if
$\mathbf{D}_2$ is nonempty and $r$ choices of colors for $\mathbf{d}_j$ if $\mathbf{D}_2$ is empty, and when $\ell\geq 2$ there are only $r-1$ choices of colors for $\mathbf{d}_j$ since $\mathbf{D}_2$ must be nonempty.

Then we establish the following recurrences
\begin{eqnarray}
P_{n+1, 1}^{(r)}  \hskip-.22cm &=&\hskip-.22cm (r-1)\sum_{k=1}^{n} P_{k, 1}^{(r)}S_{n-k}^{(r)}+(r-1)\sum_{k=0}^{n}S_{k}^{(r)}P_{n-k, 0}^{(r)}+S_{n}^{(r)}+P_{n, 1}^{(r)},   \label{eqn 3.2}  \\
P_{n+1, \ell}^{(r)}  \hskip-.22cm &=&\hskip-.22cm (r-1)\sum_{k=\ell}^{n} P_{k, \ell}^{(r)}S_{n-k}^{(r)}+(r-1)\sum_{k=0}^{n-\ell+1}S_{k}^{(r)}P_{n-k, \ell-1}^{(r)}+P_{n, \ell}^{(r)}, \  \ell\geq 2.   \label{eqn 3.3}
\end{eqnarray}

Let $P_{\ell}^{(r)}(x)=\sum_{n\geq \ell}P_{n, \ell}^{(r)}x^{n}$, by (\ref{eqn 3.1})-(\ref{eqn 3.3}), one has
\begin{eqnarray*}
P_{0}^{(r)}(x)-1  \hskip-.22cm &=&\hskip-.22cm (r-1)xP_{0}^{(r)}(x)S_r(x)+\big(S_{r}(x)-1\big)+xP_{0}^{(r)}(x), \\
P_{1}^{(r)}(x)-rx  \hskip-.22cm &=&\hskip-.22cm (r-1)xP_{1}^{(r)}(x)S_r(x)+(r-1)x\big(P_{0}^{(r)}(x)S_r(x)-1\big)  \\
\hskip-.22cm & &\hskip-.22cm \hskip1cm +x\big(S_{r}(x)-1\big)+xP_{1}^{(r)}(x), \\
P_{\ell}^{(r)}(x)-r(r-1)^{\ell-1}x^{\ell}  \hskip-.22cm &=&\hskip-.22cm (r-1)xP_{\ell}^{(r)}(x)S_r(x)+xP_{\ell}^{(r)}(x)  \\
\hskip-.22cm & &\hskip-.22cm \hskip1cm +(r-1)x\big(P_{\ell-1}^{(r)}(x)S_r(x)-r(r-1)^{\ell-2}x^{\ell-1}\big),  \  \ell\geq 2.
\end{eqnarray*}
By the relation $S_r(x)=\frac{1}{1-x-x(r-1)S_r(x)}$, after simplification, these reduce to
\begin{eqnarray*}
P_{0}^{(r)}(x)    \hskip-.22cm &=&\hskip-.22cm S_{r}(x)^2, \\
P_{1}^{(r)}(x)    \hskip-.22cm &=&\hskip-.22cm x(r-1)S_{r}(x)^2\big(\frac{1}{r-1}+S_{r}(x)^2\big), \\
P_{\ell}^{(r)}(x) \hskip-.22cm &=&\hskip-.22cm \big(x(r-1)S_{r}(x)^2\big)^{\ell}\big(\frac{1}{r-1}+S_{r}(x)^2\big),  \  \ell\geq 2,
\end{eqnarray*}
which indicate that $P_{n, \ell}^{(r)}+\frac{1}{r-1}\delta_{n, 0}$ is the $(n,\ell)$-entry of the Riordan array
$$\Big(S_r(x)^2+\frac{1}{r-1},\ x(r-1)S_r(x)^2\Big). $$
Utilizing the Lagrange inversion formula to $S_r(x)=1+xS_r(x)+x(r-1)S_r(x)^2$, that is to set $\alpha=S_r(x)-1=x(\alpha+1)((r-1)(\alpha+1)+1)$, for $n\geq 1$ and $n\geq \ell\geq 0$, one has
\begin{eqnarray*}
P_{n, \ell}^{(r)}
\hskip-.22cm &=&\hskip-.22cm [x^n]\big(x(r-1)S_{r}(x)^2\big)^{\ell}\big(\frac{1}{r-1}+S_{r}(x)^2\big) \\
\hskip-.22cm &=&\hskip-.22cm [x^{n-\ell}](r-1)^{\ell}S_{r}(x)^{2\ell+2}+[x^{n-\ell}](r-1)^{\ell-1}S_{r}(x)^{2\ell}   \\
\hskip-.22cm &=&\hskip-.22cm \frac{2\ell+2}{n-\ell}(r-1)^{\ell}[\alpha^{n-\ell-1}](\alpha+1)^{2\ell+1}\big((\alpha+1)((r-1)(\alpha+1)+1)\big)^{n-\ell} \\
\hskip-.22cm & &\hskip-.22cm  \hskip1cm +\frac{2\ell}{n-\ell}(r-1)^{\ell-1}[\alpha^{n-\ell-1}](\alpha+1)^{2\ell-1}\big((\alpha+1)((r-1)(\alpha+1)+1)\big)^{n-\ell} \\
\hskip-.22cm &=&\hskip-.22cm \sum_{j=0}^{n-\ell}\binom{n-\ell}{j}\left\{\frac{2\ell}{(2n-j)(r-1)}\binom{2n-j}{n-\ell}+\frac{2\ell+2}{2n-j+2}\binom{2n-j+2}{n-\ell} \right\}(r-1)^{n-j},
\end{eqnarray*}
as desired. \qed

\begin{corollary}\label{coro 3.1.2}
For any integers $n\geq 0$ and $r\geq 2$, we deduce
\begin{eqnarray*}
\sum_{\ell=0}^{n}(-1)^{\ell}P_{n, \ell}^{(r)}  \hskip-.22cm &=&\hskip-.22cm  S_{n}^{(r)}.
\end{eqnarray*}
\end{corollary}

\pf Note that by $S_r(x)(1-x)=1+x(r-1)S_r(x)^2$, one has
\begin{eqnarray*}
\lefteqn{\big(S_r(x)^2+\frac{1}{r-1},\ x(r-1)S_r(x)^2\big)\frac{1}{1+x}}\\
\hskip-.22cm &=&\hskip-.22cm  \big(S_r(x)^2+\frac{1}{r-1}\big)\frac{1}{1+x(r-1)S_r(x)^2}  \\
\hskip-.22cm &=&\hskip-.22cm \frac{1}{r-1}+\frac{S_r(x)^2(1-x)}{1+x(r-1)S_r(x)^2}=\frac{1}{r-1}+S_r(x),
\end{eqnarray*}
which implies that
\begin{eqnarray*}
\sum_{\ell=0}^{n}(-1)^{\ell}P_{n, \ell}^{(r)}  \hskip-.22cm &=&\hskip-.22cm [x^n]\left\{\big(S_r(x)^2+\frac{1}{r-1},\ x(r-1)S_r(x)^2\big)\frac{1}{1+x}-\frac{1}{r-1}\right\} \\
 \hskip-.22cm &=&\hskip-.22cm [x^n]S_r(x)=S_{n}^{(r)},
\end{eqnarray*}
as desired. \qed

\begin{corollary}\label{coro 3.1.3}
For any integers $n\geq 0$ and $r\geq 2$, we deduce
\begin{eqnarray*}
\sum_{\ell=1}^{n+1}(-1)^{\ell-1}\ell P_{n+1, \ell}^{(r)}  \hskip-.22cm &=&\hskip-.22cm  \sum_{\ell=0}^{n}S_{\ell+1}^{(r)}.
\end{eqnarray*}
\end{corollary}

\pf Let $P(x,y)=\sum_{n\geq 0}\big(\sum_{\ell=0}^{n}P_{n,\ell}^{(r)}y^{\ell}\big)x^{n}=1+\sum_{n\geq 0}\big(\sum_{\ell=0}^{n+1}P_{n+1,\ell}^{(r)}y^{\ell}\big)x^{n+1}$, by Theorem 3.1, one has
\begin{eqnarray*}
P(x, y) \hskip-.22cm &=&\hskip-.22cm  \big(S_r(x)^2+\frac{1}{r-1},\ x(r-1)S_r(x)^2\big)\frac{1}{1-xy}-\frac{1}{r-1} \\
\hskip-.22cm &=&\hskip-.22cm  \frac{S_r(x)^2+\frac{1}{r-1}}{1-xy(r-1)S_r(x)^2}-\frac{1}{r-1}.
\end{eqnarray*}

Note that, by $S_r(x)(1-x)=1+x(r-1)S_r(x)^2$, we deduce
\begin{eqnarray*}
\frac{\partial P(x, y)}{\partial y}\Big|_{y=-1}\hskip-.22cm &=&\hskip-.22cm   \frac{\big(S_r(x)^2+\frac{1}{r-1}\big)x(r-1)S_r(x)^2}{(1+x(r-1)S_r(x)^2)^2} \\
\hskip-.22cm &=&\hskip-.22cm  \frac{x\big((r-1)S_r(x)^2+1\big)}{(1-x)^2}=\frac{S_r(x)-1}{1-x}.
\end{eqnarray*}

Hence, we have
\begin{eqnarray*}
\sum_{\ell=1}^{n+1}(-1)^{\ell-1}\ell P_{n+1, \ell}^{(r)}  \hskip-.22cm &=&\hskip-.22cm  [x^{n+1}]\frac{\partial P(x, y)}{\partial y}\Big|_{y=-1}= [x^{n+1}]\frac{S_r(x)-1}{1-x} =\sum_{\ell=0}^{n}S_{\ell+1}^{(r)},
\end{eqnarray*}
as desired.   \qed

\section{ The statistic ``number of $\mathbf{u}$-steps" at level $\ell+1$ on $\mathcal{A}_{n+1,0}^{(r)}$ }
Let $U_{n, \ell}^{(r)}$ denote the total number of $\mathbf{u}$-steps at level $\ell+1$ in all $r$-colored Dyck paths in $\mathcal{A}_{n+1,0}^{(r)}$. The first values of $U_{n, \ell}^{(r)}$ for $r=2, 3$ are illustrated in Table 4.1 and Table 4.2 respectively.

\begin{center}
\begin{eqnarray*}
\begin{array}{c|ccccccc}\hline
n/\ell & 0      & 1      & 2       & 3       & 4    & 5    & 6    \\\hline
  0 & 2      &        &         &         &      &      &     \\
  1 & 10     & 2      &         &         &      &      &     \\
  2 & 46     & 18     & 2       &         &      &      &      \\
  3 & 214    & 118    & 26      & 2       &      &      &      \\
  4 & 1018   & 694    & 222     & 34      &  2   &      &      \\
  5 & 4946   & 3998   & 1590    & 358     &  42  &  2   &      \\
  6 & 24470  & 21434  & 10394   & 3030    &  526 &  50  & 2     \\\hline
\end{array}
\end{eqnarray*}
Table 4.1. The first values of $U_{n, \ell}^{(2)}$.
\end{center}

\begin{center}
\begin{eqnarray*}
\begin{array}{c|ccccccc}\hline
n/\ell & 0      & 1      & 2       & 3       & 4    & 5    & 6    \\\hline
  0 & 3      &        &         &         &       &      &        \\
  1 & 24     & 6      &         &         &       &      &        \\
  2 & 183    & 84     & 12      &         &       &      &        \\
  3 & 1428   & 888    & 240     & 24      &       &      &        \\
  4 & 11451  & 8580   & 3252    & 624     & 48    &      &         \\
  5 & 94020  & 79998  & 37680   & 10320   & 1536  & 96   &         \\
  6 & 787485 & 734808 & 403464  & 139728  & 30000 & 3648 & 192    \\\hline
\end{array}
\end{eqnarray*}
Table 4.2. The first values of $U_{n, \ell}^{(3)}$.
\end{center}

\begin{theorem}\label{theom 4.1.1}
For any integers $n\geq \ell\geq 0$ and $r\geq 2$, we deduce
\begin{eqnarray*}
U_{n, \ell}^{(r)}=\sum_{j=0}^{n-\ell}\binom{n-\ell}{j}\left\{\frac{2\ell+2}{2n-j+2}\binom{2n-j+2}{n-\ell}+\frac{(2\ell+3)(r-1)}{2n-j+3}\binom{2n-j+3}{n-\ell} \right\}(r-1)^{n-j}
\end{eqnarray*}
with $U_{0, 0}^{(r)}=r$. Moreover, $U_{n, \ell}^{(r)}$ is the $(n,\ell)$-entry of the Riordan array
$$\Big(S_r(x)^2\big(1+(r-1)S_r(x)\big),\ x(r-1)S_r(x)^2\Big). $$
\end{theorem}

\pf It is clear that $U_{\ell, \ell}^{(r)}=r(r-1)^{\ell}$ for $\ell\geq 0$. Given any $\mathbf{D}=\mathbf{D}_1\mathbf{u}\mathbf{D}_2\mathbf{d}_j\in \mathcal{A}_{n+2,0}^{(r)}$ with $j\in [r]$,
$\mathbf{D}_1\in \mathcal{A}_{k,0}^{(r)}$ and $\mathbf{D}_2\in \mathcal{A}_{n-k+1,0}^{(r)}$ for some $0\leq k\leq n+1$, there are two cases to be considered for the number $U_{n+1, \ell}^{(r)}$ of $\mathbf{u}$-steps at level $\ell+1$ in all $\mathbf{D}\in \mathcal{A}_{n+2,0}^{(r)}$.

When $\ell=0$, if $\mathbf{D}_1$ is empty, i.e., $k=0$, then the first $\mathbf{u}$-step of $\mathbf{D}=\mathbf{u}\mathbf{D}_2\mathbf{d}_j$ should be counted and $\mathbf{D}_2\in \mathcal{A}_{n+1,0}^{(r)}$ has $S_{n+1}^{(r)}$ choices and
$\mathbf{d}_j$ has $r-1$ choices of colors because $\mathbf{D}_2$ is nonempty. If $\mathbf{D}_1$ is nonempty, there are $U_{k-1, 0}^{(r)}$ choices of $\mathbf{u}$-steps at level $1$ in all $\mathbf{D}_1\in \mathcal{A}_{k,0}^{(r)}$ for $1\leq k\leq n+1$, and there are $S_{n-k+1}^{(r)}$ choices for $\mathbf{D}_2\in \mathcal{A}_{n-k+1,0}^{(r)}$ and $r-1$ choices of colors for $\mathbf{d}_j$ if
$\mathbf{D}_2$ is nonempty and $r$ choices of colors for $\mathbf{d}_j$ if $\mathbf{D}_2$ is empty. Then we derive the following recurrence
\begin{eqnarray}\label{eqn 4.1}
U_{n+1, 0}^{(r)}  \hskip-.22cm &=&\hskip-.22cm (r-1)\sum_{k=1}^{n+1} U_{k-1, 0}^{(r)}S_{n-k+1}^{(r)}+(r-1)\sum_{k=1}^{n+1} S_{k}^{(r)}S_{n-k+1}^{(r)}+rS_{n+1}^{(r)}+U_{n, 0}^{(r)} \nonumber  \\
                  \hskip-.22cm &=&\hskip-.22cm (r-1)\sum_{k=0}^{n} U_{k, 0}^{(r)}S_{n-k}^{(r)}+(r-1)\sum_{k=0}^{n} S_{k+1}^{(r)}S_{n-k}^{(r)}+rS_{n+1}^{(r)}+U_{n, 0}^{(r)}.
 \end{eqnarray}

When $\ell\geq 1$, we first count $\mathbf{u}$-steps at level $\ell+1$ in all $\mathbf{D}_1\in \mathcal{A}_{k,0}^{(r)}$ for $\ell+1\leq k\leq n+1$, all $\mathbf{D}_1\in \mathcal{A}_{k,0}^{(r)}$ have total $U_{k-1, \ell}^{(r)}$ $\mathbf{u}$-steps at level $\ell+1$, all $\mathbf{D}_2\in \mathcal{A}_{n-k+1,0}^{(r)}$ have $S_{n-k+1}^{(r)}$ choices and there are $r-1$ choices of colors for $\mathbf{d}_j$ if
$\mathbf{D}_2$ is nonempty and $r$ choices of colors for $\mathbf{d}_j$ if $\mathbf{D}_2$ is empty. Secondly, for the $\mathbf{u}$-steps at level $\ell$ in all $\mathbf{D}_2\in \mathcal{A}_{n-k+1,0}^{(r)}$, there are total $U_{n-k, \ell-1}^{(r)}$ $\mathbf{u}$-steps at level $\ell$ for $0\leq k\leq n-\ell+1$, and all $\mathbf{D}_1\in \mathcal{A}_{k,0}^{(r)}$ have $S_{k}^{(r)}$ choices. Moreover, there are $r-1$ choices of colors for $\mathbf{d}_j$ since $\mathbf{D}_2$ must be nonempty.

Then we establish the following recurrence
\begin{eqnarray}\label{eqn 4.2}
U_{n+1, \ell}^{(r)}  \hskip-.22cm &=&\hskip-.22cm (r-1)\sum_{k=\ell+1}^{n+1} U_{k-1, \ell}^{(r)}S_{n-k+1}^{(r)}+(r-1)\sum_{k=0}^{n-\ell+1}S_{k}^{(r)}U_{n-k, \ell-1}^{(r)}+U_{n, \ell}^{(r)},  \nonumber \\
\hskip-.22cm &=&\hskip-.22cm (r-1)\sum_{k=\ell}^{n} U_{k, \ell}^{(r)}S_{n-k}^{(r)}+(r-1)\sum_{k=0}^{n-\ell+1}S_{k}^{(r)}U_{n-k, \ell-1}^{(r)}+U_{n, \ell}^{(r)}, \hskip0.8cm  \ell\geq 1.
\end{eqnarray}

Let $U_{\ell}^{(r)}(x)=\sum_{n\geq \ell}U_{n, \ell}^{(r)}x^{n}$, by (\ref{eqn 4.1})-(\ref{eqn 4.2}), one has
\begin{eqnarray*}
U_{0}^{(r)}(x)-r  \hskip-.22cm &=&\hskip-.22cm (r-1)xU_{0}^{(r)}(x)S_r(x)+(r-1)\big(S_{r}(x)-1\big)S_r(x) \\
\hskip-.22cm & &\hskip-.22cm        \hskip0.5cm    +r\big(S_{r}(x)-1\big)+xU_{0}^{(r)}(x), \\
U_{\ell}^{(r)}(x)-r(r-1)^{\ell}x^{\ell}  \hskip-.22cm &=&\hskip-.22cm (r-1)xU_{\ell}^{(r)}(x)S_r(x)+xU_{\ell}^{(r)}(x)  \\
\hskip-.22cm & &\hskip-.22cm \hskip0.5cm +(r-1)x\big(U_{\ell-1}^{(r)}(x)S_r(x)-r(r-1)^{\ell-1}x^{\ell-1}\big),  \  \ell\geq 1.
\end{eqnarray*}
By the relation $S_r(x)=\frac{1}{1-x-x(r-1)S_r(x)}$, after simplification, these generate
\begin{eqnarray*}
U_{0}^{(r)}(x)    \hskip-.22cm &=&\hskip-.22cm S_{r}(x)^2\big(1+(r-1)S_r(x)\big), \\
U_{\ell}^{(r)}(x) \hskip-.22cm &=&\hskip-.22cm \big(x(r-1)S_{r}(x)^2\big)^{\ell}S_{r}(x)^2\big(1+(r-1)S_r(x)\big),  \  \ell\geq 1,
\end{eqnarray*}
which indicate that $U_{n, \ell}^{(r)}$ is the $(n,\ell)$-entry of the Riordan array
$$\Big(S_{r}(x)^2\big(1+(r-1)S_r(x)\big),\ x(r-1)S_r(x)^2\Big). $$
By using the Lagrange inversion formula to $S_r(x)=1+xS_r(x)+x(r-1)S_r(x)^2$, that is to set $\alpha=S_r(x)-1=x(\alpha+1)((r-1)(\alpha+1)+1)$, for $n\geq \ell\geq 0$, one has
\begin{eqnarray*}
U_{n, \ell}^{(r)}
\hskip-.22cm &=&\hskip-.22cm [x^n]\big(x(r-1)S_{r}(x)^2\big)^{\ell} S_{r}(x)^2\big(1+(r-1)S_r(x)\big) \\
\hskip-.22cm &=&\hskip-.22cm [x^{n-\ell}](r-1)^{\ell}S_{r}(x)^{2\ell+2}+[x^{n-\ell}](r-1)^{\ell+1}S_{r}(x)^{2\ell+3}   \\
\hskip-.22cm &=&\hskip-.22cm \frac{2\ell+2}{n-\ell}(r-1)^{\ell}[\alpha^{n-\ell-1}](\alpha+1)^{2\ell+1}\big((\alpha+1)((r-1)(\alpha+1)+1)\big)^{n-\ell} \\
\hskip-.22cm & &\hskip-.22cm  \hskip1cm +\frac{2\ell+3}{n-\ell}(r-1)^{\ell+1}[\alpha^{n-\ell-1}](\alpha+1)^{2\ell+2}\big((\alpha+1)((r-1)(\alpha+1)+1)\big)^{n-\ell} \\
\hskip-.22cm &=&\hskip-.22cm \sum_{j=0}^{n-\ell}\binom{n-\ell}{j}\left\{\frac{2\ell+2}{2n-j+2}\binom{2n-j+2}{n-\ell}+\frac{(2\ell+3)(r-1)}{2n-j+3}\binom{2n-j+3}{n-\ell} \right\}(r-1)^{n-j},
\end{eqnarray*}
as desired. \qed

\begin{corollary}\label{coro 4.1.2}
For any integers $n\geq 0$ and $r\geq 2$, we deduce
\begin{eqnarray*}
\sum_{\ell=0}^{n}(-1)^{\ell}U_{n, \ell}^{(r)}S_{\ell}^{(r)}  \hskip-.22cm &=&\hskip-.22cm  \sum_{\ell=0}^{n}(-1)^{n-\ell}S_{\ell+1}^{(r)}S_{n-\ell}(-1, (r-1)^2).
\end{eqnarray*}
Specially, the case $r=2$ generates that
\begin{eqnarray*}
\sum_{\ell=0}^{n}(-1)^{\ell}U_{n, \ell}^{(2)}S_{\ell}  \hskip-.22cm &=&\hskip-.22cm  S_{n+1}.
\end{eqnarray*}
\end{corollary}

\pf Note that by $\frac{S_{r}(x)-1}{x}=S_r(x)(1+(r-1)S_r(x))$ and (\ref{eqn 2.5}), one has
\begin{eqnarray*}
\lefteqn{\big(S_{r}(x)^2\big(1+(r-1)S_r(x)\big),\ x(r-1)S_r(x)^2\big)S_r(-x) }\\
\hskip-.22cm &=&\hskip-.22cm  S_{r}(x)^2\big(1+(r-1)S_r(x)\big)S_r(-x(r-1)S_r(x)^2)  \\
\hskip-.22cm &=&\hskip-.22cm  S_{r}(x)\big(1+(r-1)S_r(x)\big)S(-1, (r-1)^2; -x) \\
\hskip-.22cm &=&\hskip-.22cm \frac{S_{r}(x)-1}{x}S(-1, (r-1)^2; -x),
\end{eqnarray*}
which implies that
\begin{eqnarray*}
\sum_{\ell=0}^{n}(-1)^{\ell}U_{n, \ell}^{(r)}S_{\ell}^{(r)}  \hskip-.22cm &=&\hskip-.22cm [x^n]\Big(S_{r}(x)^2\big(1+(r-1)S_r(x)\big),\ x(r-1)S_r(x)^2\Big)S_r(-x) \\
 \hskip-.22cm &=&\hskip-.22cm [x^n]\frac{S_{r}(x)-1}{x}S(-1, (r-1)^2; -x) \\
 \hskip-.22cm &=&\hskip-.22cm \sum_{\ell=0}^{n}(-1)^{n-\ell}S_{\ell+1}^{(r)}S_{n-\ell}(-1, (r-1)^2),
\end{eqnarray*}
as desired.

Specially, in the case $r=2$, by $S_{\ell}^{(2)}=S_{\ell}$ and $S_{n-\ell}(-1, 1)=\delta_{n-\ell, 0}$, one has
\begin{eqnarray*}
\sum_{\ell=0}^{n}(-1)^{\ell}U_{n, \ell}^{(2)}S_{\ell}  \hskip-.22cm &=&\hskip-.22cm  S_{n+1}.
\end{eqnarray*}
This completes the proof.  \qed

\begin{corollary}\label{coro 4.1.3}
For any integers $n\geq 0$ and $r\geq 2$, we deduce
\begin{eqnarray*}
\lefteqn{\sum_{\ell=0}^{n}(-1)^{\ell}U_{n, \ell}^{(r)}\big(S_{\ell+1}^{(r)}-S_{\ell}^{(r)}\big) } \\
\hskip-.22cm &=&\hskip-.22cm  \sum_{\ell=0}^{n}(-1)^{n-\ell}\Big(\frac{\delta_{\ell, 0}}{r-1}+S_{\ell}^{(r)}\Big)\big(S_{n-\ell}(-1, (r-1)^2)+S_{n-\ell+1}(-1, (r-1)^2)\big).
\end{eqnarray*}
Specially, the case $r=2$ generates that
\begin{eqnarray*}
\sum_{\ell=0}^{n}(-1)^{\ell}U_{n, \ell}^{(2)}\big(S_{\ell+1}-S_{\ell}\big)  \hskip-.22cm &=&\hskip-.22cm  \delta_{n, 0}+S_{n}.
\end{eqnarray*}
\end{corollary}

\pf Note that by the relation
$$S_r(x)^2=\frac{S_r(x)-xS_r(x)-1}{x(r-1)}=\sum_{n\geq 0}\frac{S_{n+1}^{(r)}-S_{n}^{(r )}}{r-1}x^n$$
and (\ref{eqn 2.5}), and by (\ref{eqn 1.4}) in the case $a=-1, b=(r-1)^2$, one has
\begin{eqnarray*}
\lefteqn{\big(S_{r}(x)^2\big(1+(r-1)S_r(x)\big),\ x(r-1)S_r(x)^2\big)S_r(-x)^2 }\\
\hskip-.22cm &=&\hskip-.22cm  S_{r}(x)^2\big(1+(r-1)S_r(x)\big)S_r(-x(r-1)S_r(x)^2)^2  \\
\hskip-.22cm &=&\hskip-.22cm  \big(1+(r-1)S_r(x)\big)S(-1, (r-1)^2; -x)^2 \\
\hskip-.22cm &=&\hskip-.22cm  \big(1+(r-1)S_r(x)\big)\frac{1-S(-1, (r-1)^2; -x)-xS(-1, (r-1)^2; -x)}{x(r-1)^2} \\
\hskip-.22cm &=&\hskip-.22cm  \Big(\frac{1}{r-1}+S_r(x)\Big)\frac{1-S(-1, (r-1)^2; -x)+xS(-1, (r-1)^2; -x)}{x(r-1)},
\end{eqnarray*}
which implies that
\begin{eqnarray*}
\lefteqn{\sum_{\ell=0}^{n}(-1)^{\ell}U_{n, \ell}^{(r)}\big(S_{\ell+1}^{(r)}-S_{\ell}^{(r)}\big)} \\
 \hskip-.22cm &=&\hskip-.22cm [x^n](r-1)\Big(S_{r}(x)^2\big(1+(r-1)S_r(x)\big),\ x(r-1)S_r(x)^2\Big)S_r(-x)^2 \\
 \hskip-.22cm &=&\hskip-.22cm [x^n](r-1)\Big(\frac{1}{r-1}+S_r(x)\Big)\frac{1-S(-1, (r-1)^2; -x)+xS(-1, (r-1)^2; -x)}{x(r-1)}  \\
 \hskip-.22cm &=&\hskip-.22cm \sum_{\ell=0}^{n}(-1)^{n-\ell}\Big(\frac{\delta_{\ell, 0}}{r-1}+S_{\ell}^{(r)}\Big)\big(S_{n-\ell}(-1, (r-1)^2)+S_{n-\ell+1}(-1, (r-1)^2)\big),
\end{eqnarray*}
as desired.

Specially, in the case $r=2$, by $S_{\ell}^{(2)}=S_{\ell}$ and $S_{n-\ell}(-1, 1)=\delta_{n-\ell, 0}$, one has
\begin{eqnarray*}
\sum_{\ell=0}^{n}(-1)^{\ell}U_{n, \ell}^{(2)}\big(S_{\ell+1}-S_{\ell}\big)  \hskip-.22cm &=&\hskip-.22cm  \delta_{n, 0}+S_{n}.
\end{eqnarray*}
This completes the proof.  \qed

\begin{corollary}\label{coro 4.1.4}
For any integers $n\geq 0$ and $r\geq 2$, we deduce
\begin{eqnarray*}
\sum_{\ell=0}^{n}(-1)^{\ell}(\ell+1)U_{n, \ell}^{(r)}  \hskip-.22cm &=&\hskip-.22cm (n+1)+(r-1)\sum_{\ell=0}^{n}(\ell+1)S_{n-\ell}^{(r)}.
\end{eqnarray*}
\end{corollary}

\pf Note that by $S_r(x)(1-x)=1+x(r-1)S_r(x)^2$, one has
\begin{eqnarray*}
\lefteqn{\big(S_{r}(x)^2\big(1+(r-1)S_r(x)\big),\ x(r-1)S_r(x)^2\big)\frac{1}{(1+x)^2}}\\
\hskip-.22cm &=&\hskip-.22cm  S_{r}(x)^2\big(1+(r-1)S_r(x)\big)\frac{1}{(1+x(r-1)S_r(x)^2)^2}  \\
\hskip-.22cm &=&\hskip-.22cm S_{r}(x)^2\big(1+(r-1)S_r(x)\big)\frac{1}{S_r(x)^2(1-x)^2} \\
\hskip-.22cm &=&\hskip-.22cm \frac{1+(r-1)S_r(x)}{(1-x)^2},
\end{eqnarray*}
which implies that
\begin{eqnarray*}
\sum_{\ell=0}^{n}(-1)^{\ell}(\ell+1)U_{n, \ell}^{(r)}  \hskip-.22cm &=&\hskip-.22cm [x^n]\left\{\big(S_{r}(x)^2\big(1+(r-1)S_r(x)\big),\ x(r-1)S_r(x)^2\big)\frac{1}{(1+x)^2} \right\}\\
 \hskip-.22cm &=&\hskip-.22cm [x^n]\frac{1+(r-1)S_r(x)}{(1-x)^2}  \\
 \hskip-.22cm &=&\hskip-.22cm (n+1)+(r-1)\sum_{\ell=0}^{n}(\ell+1)S_{n-\ell}^{(r)},
\end{eqnarray*}
as desired. \qed

\section{ The statistic ``number of peaks" at level $\ell+1$ on $\mathcal{A}_{n+1,0}^{(r)}$ }
Let $p_{n, \ell}^{(r)}$ denote the total number of peaks at level $\ell+1$ in all $r$-colored Dyck paths in $\mathcal{A}_{n+1,0}^{(r)}$. The first values of $p_{n, \ell}^{(r)}$ for $r=2, 3$ are illustrated in Table 5.1 and Table 5.2 respectively.

\begin{center}
\begin{eqnarray*}
\begin{array}{c|ccccccc}\hline
n/\ell & 0      & 1      & 2       & 3       & 4    & 5    & 6    \\\hline
  0 & 2      &        &         &         &      &      &     \\
  1 & 8      & 2      &         &         &      &      &     \\
  2 & 32     & 16     & 2       &         &      &      &      \\
  3 & 136    & 96     & 24      & 2       &      &      &      \\
  4 & 608    & 528    & 192     & 32      &  2   &      &      \\
  5 & 2824   & 2816   & 1304    & 320     &  40  &  2   &      \\
  6 & 13504  & 14864  & 8160    & 2592    &  480 &  48  & 2     \\\hline
\end{array}
\end{eqnarray*}
Table 5.1. The first values of $p_{n, \ell}^{(2)}$.
\end{center}

\begin{center}
\begin{eqnarray*}
\begin{array}{c|ccccccc}\hline
n/\ell & 0      & 1      & 2       & 3       & 4    & 5    & 6    \\\hline
  0 & 3      &        &         &         &       &      &        \\
  1 & 18     & 6      &         &         &       &      &        \\
  2 & 117    & 72     & 12      &         &       &      &        \\
  3 & 828    & 684    & 216     & 24      &       &      &        \\
  4 & 6219   & 6120   & 2700    & 576     & 48    &      &         \\
  5 & 48726  & 53874  & 29376   & 8928    & 1440  & 96   &         \\
  6 & 393885 & 473328 & 299160  & 114624  & 26640 & 3456 & 192    \\\hline
\end{array}
\end{eqnarray*}
Table 5.2. The first values of $p_{n, \ell}^{(3)}$.
\end{center}

\begin{theorem}\label{theom 5.1.1}
For any integers $n\geq \ell\geq 0$ and $r\geq 2$, we deduce
\begin{eqnarray*}
p_{n, \ell}^{(r)}=\sum_{j=0}^{n-\ell}\frac{2r(\ell+1)}{2n-j+2}\binom{2n-j+2}{n-\ell}\binom{n-\ell}{j}(r-1)^{n-j}
\end{eqnarray*}
with $p_{0, 0}^{(r)}=r$. Moreover, $p_{n, \ell}^{(r)}$ is the $(n,\ell)$-entry of the Riordan array
$$\Big(rS_r(x)^2,\ x(r-1)S_r(x)^2\Big). $$
\end{theorem}

\pf Clearly, we have $p_{\ell, \ell}^{(r)}=r(r-1)^{\ell}$ for $\ell\geq 0$. Given any $\mathbf{D}=\mathbf{D}_1\mathbf{u}\mathbf{D}_2\mathbf{d}_j\in \mathcal{A}_{n+2,0}^{(r)}$ with $j\in [r]$,
$\mathbf{D}_1\in \mathcal{A}_{k,0}^{(r)}$ and $\mathbf{D}_2\in \mathcal{A}_{n-k+1,0}^{(r)}$ for some $0\leq k\leq n+1$, there are two cases to be considered for the number $p_{n+1, \ell}^{(r)}$ of peaks at level $\ell+1$ in all $r$-colored Dyck paths $\mathbf{D}\in \mathcal{A}_{n+2,0}^{(r)}$.

When $\ell=0$, we first count peaks at level $1$ in all $r$-colored Dyck paths $\mathbf{D}_1\in \mathcal{A}_{k,0}^{(r)}$ for $1\leq k\leq n+1$, since $\mathbf{D}_1$ is nonempty, there are $p_{k-1, 0}^{(r)}$ choices of peaks at level $1$ in all $\mathbf{D}_1\in \mathcal{A}_{k,0}^{(r)}$ for $1\leq k\leq n+1$, and there are $S_{n-k+1}^{(r)}$ choices for $\mathbf{D}_2\in \mathcal{A}_{n-k+1,0}^{(r)}$ and $r-1$ choices of colors for $\mathbf{d}_j$ if
$\mathbf{D}_2$ is nonempty and $r$ choices of colors for $\mathbf{d}_j$ if $\mathbf{D}_2$ is empty. Secondly, for the peaks at level $1$ in all $\mathbf{u}\mathbf{D}_2\mathbf{d}_j$, which only occurs when $\mathbf{D}_2$ is empty, so there are $r$ choices of colors for $\mathbf{d}_j$ and $S_{n+1}^{(r)}$ choices for $\mathbf{D}_1\in \mathcal{A}_{n+1,0}^{(r)}$. Then we obtain the following recurrence
\begin{eqnarray}\label{eqn 5.1}
p_{n+1, 0}^{(r)}  \hskip-.22cm &=&\hskip-.22cm (r-1)\sum_{k=1}^{n+1} p_{k-1, 0}^{(r)}S_{n-k+1}^{(r)}+p_{n, 0}^{(r)}+rS_{n+1}^{(r)}   \nonumber \\
\hskip-.22cm &=&\hskip-.22cm (r-1)\sum_{k=0}^{n} p_{k, 0}^{(r)}S_{n-k}^{(r)}+p_{n, 0}^{(r)}+rS_{n+1}^{(r)}.
 \end{eqnarray}

When $\ell\geq 1$, let's first count peaks at level $\ell+1$ in all $r$-colored Dyck paths $\mathbf{D}_1\in \mathcal{A}_{k,0}^{(r)}$ for $\ell+1\leq k\leq n+1$, all $\mathbf{D}_1\in \mathcal{A}_{k,0}^{(r)}$ have total $p_{k-1, \ell}^{(r)}$ peaks at level $\ell+1$, all $\mathbf{D}_2\in \mathcal{A}_{n-k+1,0}^{(r)}$ have $S_{n-k+1}^{(r)}$ choices and there are $r-1$ choices of colors for $\mathbf{d}_j$ if
$\mathbf{D}_2$ is nonempty and $r$ choices of colors for $\mathbf{d}_j$ if $\mathbf{D}_2$ is empty. Secondly, for the peaks at level $\ell$ in all $r$-colored Dyck paths $\mathbf{D}_2\in \mathcal{A}_{n-k+1,0}^{(r)}$, there are total $p_{n-k, \ell-1}^{(r)}$ peaks at level $\ell$ for $0\leq k\leq n-\ell+1$, and all $\mathbf{D}_1\in \mathcal{A}_{k,0}^{(r)}$ have $S_{k}^{(r)}$ choices. Moreover, there are $r-1$ choices of colors for $\mathbf{d}_j$ since $\mathbf{D}_2$ must be nonempty.

Hence, we deduce the following recurrence
\begin{eqnarray}\label{eqn 5.2}
p_{n+1, \ell}^{(r)}  \hskip-.22cm &=&\hskip-.22cm  (r-1)\sum_{k=\ell+1}^{n+1} p_{k-1, \ell}^{(r)}S_{n-k+1}^{(r)}+(r-1)\sum_{k=0}^{n-\ell+1}S_{k}^{(r)}p_{n-k, \ell-1}^{(r)}+p_{n, \ell}^{(r)}   \nonumber  \\
\hskip-.22cm &=&\hskip-.22cm  (r-1)\sum_{k=\ell}^{n} p_{k, \ell}^{(r)}S_{n-k}^{(r)}+(r-1)\sum_{k=0}^{n-\ell+1}S_{k}^{(r)}p_{n-k, \ell-1}^{(r)}+p_{n, \ell}^{(r)}, \   \ \ \  \ell\geq 1.
\end{eqnarray}

Let $p_{\ell}^{(r)}(x)=\sum_{n\geq \ell}p_{n, \ell}^{(r)}x^{n}$, by (\ref{eqn 5.1})-(\ref{eqn 5.2}), one has
\begin{eqnarray*}
p_{0}^{(r)}(x)-r  \hskip-.22cm &=&\hskip-.22cm (r-1)xp_{0}^{(r)}(x)S_r(x)+xp_{0}^{(r)}(x)+r\big(S_{r}(x)-1\big), \\
p_{\ell}^{(r)}(x)-r(r-1)^{\ell}x^{\ell}  \hskip-.22cm &=&\hskip-.22cm (r-1)xp_{\ell}^{(r)}(x)S_r(x)+xp_{\ell}^{(r)}(x)  \\
\hskip-.22cm & &\hskip-.22cm \hskip0.5cm +(r-1)x\big(p_{\ell-1}^{(r)}(x)S_r(x)-r(r-1)^{\ell-1}x^{\ell-1}\big),  \  \ell\geq 1.
\end{eqnarray*}
By the relation $S_r(x)=\frac{1}{1-x-x(r-1)S_r(x)}$, after simplification, these generate
\begin{eqnarray*}
p_{0}^{(r)}(x)    \hskip-.22cm &=&\hskip-.22cm rS_{r}(x)^2, \\
p_{\ell}^{(r)}(x) \hskip-.22cm &=&\hskip-.22cm rS_{r}(x)^2\big(x(r-1)S_{r}(x)^2\big)^{\ell},  \  \ell\geq 1,
\end{eqnarray*}
which indicate that $p_{n, \ell}^{(r)}$ is the $(n,\ell)$-entry of the Riordan array
$$\Big(rS_{r}(x)^2,\ x(r-1)S_r(x)^2\Big). $$
By using the Lagrange inversion formula to $S_r(x)=1+xS_r(x)+x(r-1)S_r(x)^2$, one has
\begin{eqnarray*}
p_{n, \ell}^{(r)} \hskip-.22cm &=&\hskip-.22cm   [x^n] rS_{r}(x)^2\big(x(r-1)S_{r}(x)^2\big)^{\ell}  \\
                  \hskip-.22cm &=&\hskip-.22cm   \sum_{j=0}^{n-\ell}\frac{2r(\ell+1)}{2n-j+2}\binom{2n-j+2}{n-\ell}\binom{n-\ell}{j}(r-1)^{n-j},
\end{eqnarray*}
as desired. The detail is omitted. \qed

Similar to the proof of Corollary \ref{coro 4.1.2} or \ref{coro 4.1.3}, we obtain the following result, the detail is left to the interested readers.
\begin{corollary}\label{coro 5.1.2}
For any integers $n\geq 0$ and $r\geq 2$, we deduce
\begin{eqnarray*}
\sum_{\ell=0}^{n}(-1)^{\ell}p_{n, \ell}^{(r)}S_{\ell}^{(r)}  \hskip-.22cm &=&\hskip-.22cm  r\sum_{\ell=0}^{n}(-1)^{n-\ell}S_{\ell}^{(r)}S_{n-\ell}(-1, (r-1)^2), \\
\sum_{\ell=0}^{n}(-1)^{\ell}p_{n, \ell}^{(r)}\big(S_{\ell+1}^{(r)}-S_{\ell}^{(r)}\big)  \hskip-.22cm &=&\hskip-.22cm  \frac{r}{r-1}(-1)^{n}\big(S_{n+1}(-1, (r-1)^2)+S_{n}(-1, (r-1)^2)\big).
\end{eqnarray*}
Specially, the case $r=2$ produces respectively that
\begin{eqnarray*}
\sum_{\ell=0}^{n}(-1)^{\ell}p_{n, \ell}^{(2)}S_{\ell}  \hskip-.22cm &=&\hskip-.22cm  2S_{n}, \\
\sum_{\ell=0}^{n}(-1)^{\ell}p_{n, \ell}^{(2)}\big(S_{\ell+1}-S_{\ell}\big)  \hskip-.22cm &=&\hskip-.22cm  2\delta_{n,0}.
\end{eqnarray*}
\end{corollary}

\begin{corollary}\label{coro 5.1.3}
Let $m$ be any complex number, for any integers $n\geq 0$ and $r\geq 2$, we deduce
\begin{eqnarray*}
\lefteqn{\sum_{\ell=0}^{n}(-1)^{\ell}\binom{m+\ell}{\ell} p_{n, \ell}^{(r)}  } \\
\hskip-.22cm &=&\hskip-.22cm  \sum_{\ell=0}^{n}\sum_{j=0}^{n-\ell}\binom{m+\ell}{\ell}\binom{n-\ell}{j}\frac{(1-m)r}{n-\ell+j-m+1}\binom{n-\ell+j-m+1}{n-\ell}(r-1)^j.
\end{eqnarray*}
Specially, the cases $m=0, 1$ generate respectively that
\begin{eqnarray*}
\sum_{\ell=0}^{n}(-1)^{\ell}p_{n, \ell}^{(r)} \hskip-.22cm &=&\hskip-.22cm  r\sum_{\ell=0}^{n}S_{\ell}^{(r)},  \\
\sum_{\ell=0}^{n}(-1)^{\ell}(\ell +1)p_{n, \ell}^{(r)} \hskip-.22cm &=&\hskip-.22cm  r(n+1).
\end{eqnarray*}
\end{corollary}

\pf Note that by $S_r(x)(1-x)=1+x(r-1)S_r(x)^2$, one has
\begin{eqnarray*}
\lefteqn{\big(rS_{r}(x)^2,\ x(r-1)S_r(x)^2\big)\frac{1}{(1+x)^{m+1}}   }       \\
\hskip-.22cm &=&\hskip-.22cm  rS_{r}(x)^2\frac{1}{(1+x(r-1)S_r(x)^2)^{m+1}} = \frac{rS_r(x)^{1-m}}{(1-x)^{m+1}}.
\end{eqnarray*}
By using the Lagrange inversion formula to $S_r(x)=1+xS_r(x)+x(r-1)S_r(x)^2$, that is to set $\alpha=S_r(x)-1=x(\alpha+1)((r-1)(\alpha+1)+1)$, for $n\geq \ell\geq 0$, one has
\begin{eqnarray*}
\lefteqn{ \sum_{\ell=0}^{n}(-1)^{\ell}\binom{m+\ell}{\ell} p_{n, \ell}^{(r)}  } \\
 \hskip-.22cm &=&\hskip-.22cm [x^n]\left\{\big(rS_{r}(x)^2,\ x(r-1)S_r(x)^2\big)\frac{1}{(1+x)^{m+1}} \right\}\\
 \hskip-.22cm &=&\hskip-.22cm [x^n]\frac{rS_r(x)^{1-m}}{(1-x)^{m+1}} = \sum_{\ell=0}^{n}\binom{m+\ell}{\ell}[x^{n-\ell}]rS_r(x)^{1-m} \\
 \hskip-.22cm &=&\hskip-.22cm \sum_{\ell=0}^{n}\sum_{j=0}^{n-\ell}\binom{m+\ell}{\ell}\binom{n-\ell}{j}\frac{(1-m)r}{n-\ell+j-m+1}\binom{n-\ell+j-m+1}{n-\ell}(r-1)^j,
\end{eqnarray*}
as desired.

Specially, one could verify easily the results when $m=0$ and $m=1$.  \qed

\section{ The statistic ``number of $\mathbf{udu}$-steps" on $\mathcal{A}_{n,0}^{(r)}$  }
Let $\mathcal{T}_{n, \ell}^{(r)}$ denote the subset of all $r$-colored Dyck paths in $\mathcal{A}_{n,0}^{(r)}$ with exact $\ell$ $\mathbf{udu}$-steps and $T_{n, \ell}^{(r)}$ be its cardinality. The first values of $T_{n, \ell}^{(r)}$ for $r=2$ and $r=3$ are illustrated in Table 6.1 and 6.2 respectively.

\begin{center}
\begin{eqnarray*}
\begin{array}{c|ccccccc}\hline
n/\ell & 0   & 1      & 2       & 3       & 4    & 5       \\\hline
  0 & 1      &        &         &         &      &       \\
  1 & 2      &        &         &         &      &       \\
  2 & 2      & 4      &         &         &      &        \\
  3 & 6      & 8      & 8       &         &      &        \\
  4 & 14     & 36     & 24      & 16      &      &        \\
  5 & 42     & 112    & 144     & 64      & 32   &        \\
  6 & 122    & 420    & 560     & 480     & 160  & 64     \\\hline
\end{array}
\end{eqnarray*}
Table 6.1. The first values of $T_{n, \ell}^{(2)}$.
\end{center}

\begin{center}
\begin{eqnarray*}
\begin{array}{c|ccccccc}\hline
n/\ell & 0   & 1      & 2       & 3       & 4    & 5       \\\hline
  0 & 1      &        &         &         &      &       \\
  1 & 3      &        &         &         &      &       \\
  2 & 6      & 9      &         &         &      &        \\
  3 & 30     & 36     & 27      &         &      &        \\
  4 & 132    & 270    & 162     & 81      &      &        \\
  5 & 696    & 1584   & 1620    & 648     & 243  &        \\
  6 & 3696   & 10440  & 11880   & 8100    & 2430 &  729    \\\hline
\end{array}
\end{eqnarray*}
Table 6.2. The first values of $T_{n, \ell}^{(3)}$.
\end{center}

\begin{theorem}\label{theom 6.1.1}
For any integers $n> \ell\geq 0$ and $r\geq 2$, we deduce
\begin{eqnarray}\label{eqn 6.1}
T_{n, \ell}^{(r)}= \binom{n-1}{\ell}T_{n-\ell,0}^{(r)}r^{\ell}
\end{eqnarray}
with $T_{0, 0}^{(r)}=1$, where
\begin{eqnarray}\label{eqn 6.2}
T_{n,0}^{(r)}=\sum_{j=[\frac{n}{2}]}^{n-1}\binom{n-1}{j}\frac{1}{j+1}\binom{j+1}{n-j}r^{n-j}(r-1)^{j}, \ (n\geq 1),
\end{eqnarray}
is determinated by
\begin{eqnarray}\label{eqn 6.3}
T_0(x)=\sum_{n\geq 0}T_{n,0}^{(r)}x^{n}=\frac{1+(r-1)x-\sqrt{(1+(r-1)x)^2-4(r-1)(1+rx)x}}{2(r-1)x},
\end{eqnarray}
which satisfies the recurrence $T_0(x)=1+rx-(r-1)xT_0(x)+(r-1)xT_0(x)^2$.
\end{theorem}

\pf First we consider $\mathcal{T}_{n,0}^{(r)}$, the subset of $r$-colored Dyck paths of length $2n$ with no $\mathbf{udu}$-steps (called $\mathbf{udu}$-avoiding).
For any $\mathbf{D}\in \mathcal{T}_{n,0}^{(r)}$, $\mathbf{D}$ has one of the following three forms
$$ \varepsilon,\ \mathbf{ud}_j\ \mbox{or}\ \mathbf{u}\mathbf{D}_1\mathbf{d}_j\mathbf{D}_2,$$
where both $\mathbf{D}_1(\neq \varepsilon)$ and $\mathbf{D}_2$ are $\mathbf{udu}$-avoiding, $\mathbf{d}_j$ has color $j$ for some $j\in [r]$, so $\mathbf{d}_j$ has $r$ choices for colors in the second case and has $r-1$ choices in the third case since $\mathbf{D}_1$ is nonempty. Let $T_0(x)=\sum_{n\geq 0}T_{n,0}^{(r)}x^{n}$. Then
\begin{eqnarray*}
T_0(x)=1+rx+(r-1)x(T_0(x)-1)T_0(x).
\end{eqnarray*}
By solving it, one obtains (\ref{eqn 6.3}).

By using the Lagrange inversion formula to $T_0(x)=1+rx+(r-1)x(T_0(x)-1)T_0(x)$, that is to set $\beta=T_0(x)-1=x(r+(r-1)\beta(\beta+1))$, for $n>0$, one has
\begin{eqnarray*}
T_{n,0}^{(r)} \hskip-.22cm &=&\hskip-.22cm [x^n](T_0(x)-1)=\frac{1}{n}[\beta^{n-1}]\big(r+(r-1)\beta(\beta+1)\big)^{n}  \\
  \hskip-.22cm &=&\hskip-.22cm \sum_{j=[\frac{n}{2}]}^{n-1}\binom{n-1}{j}\frac{1}{j+1}\binom{j+1}{n-j}r^{n-j}(r-1)^{j},
\end{eqnarray*}
which leads to (\ref{eqn 6.2}).

For any $\mathbf{D}\in \mathcal{T}_{n,\ell}^{(r)}$, we can get a $\mathbf{D}^{*}\in \mathcal{T}_{n-\ell,0}^{(r)}$ by deleting the $\ell$ $\mathbf{d}_j\mathbf{u}$-steps once they form $\mathbf{udu}$-steps, where $\mathbf{d}_j$ has color $j$ for some $j\in [r]$. Conversely, for any $\mathbf{D}^{*}\in \mathcal{T}_{n-\ell,0}^{(r)}$ there are exact $n-\ell$ endpoints of $\mathbf{u}$-steps, one can repeatedly select $\ell$ endpoints and insert $\ell$ $\mathbf{d}_j\mathbf{u}$-steps to form $\ell$ $\mathbf{u}\mathbf{d}\mathbf{u}$-steps. Since there are $\binom{n-1}{\ell}$ ways to select repeatedly $\ell$ endpoints from $n-\ell$ endpoints, and each $\mathbf{d}_j$ has one of $r$ colors, then we obtain that
\begin{eqnarray*}
T_{n, \ell}^{(r)}= \binom{n-1}{\ell}T_{n-\ell,0}^{(r)}r^{\ell},
\end{eqnarray*}
which proves (\ref{eqn 6.1}).  \qed

\begin{theorem}\label{theom 6.1.2}
For any integers $n\geq 1$ and $r\geq 2$, we deduce
\begin{eqnarray*}
\sum_{\ell=1}^{n} \ell T_{n+1, \ell}^{(r)}  \hskip-.22cm &=&\hskip-.22cm  rnS_n^{(r)}.
\end{eqnarray*}
\end{theorem}

\pf Let $T_r(x,y)=1+\sum_{n\geq 1}\big(\sum_{\ell=0}^{n-1}T_{n,\ell}^{(r)}y^{\ell}\big)x^{n}=1+\sum_{n\geq 0}\big(\sum_{\ell=0}^{n}T_{n+1,\ell}^{(r)}y^{\ell}\big)x^{n+1}$, similar to obtaining the relation for $T_0(x)$ in Theorem 6.1, one has
\begin{eqnarray*}
T_r(x,y)=1+rx+rxy(T_r(x,y)-1)+(r-1)x(T_r(x,y)-1)T_r(x,y),
\end{eqnarray*}
which, when $y=0$ or $1$, generates respectively that $T_r(x,0)=T_0(x)$ and $T_r(x,1)=S_r(x)$.

Note that
\begin{eqnarray*}
\frac{\partial T_r(x, y)}{\partial y}\Big|_{y=1}=\frac{rx(T_r(x,1)-1)}{1-x-2(r-1)xT_r(x,1)}=rx^2\frac{\partial T_r(x, 1)}{\partial x}.
\end{eqnarray*}

Hence, we deduce
\begin{eqnarray*}
\sum_{\ell=1}^{n} \ell T_{n+1, \ell}^{(r)}  \hskip-.22cm &=&\hskip-.22cm  [x^{n+1}]\frac{\partial T_r(x, y)}{\partial y}\Big|_{y=1}=[x^{n+1}]rx^2\frac{\partial T_r(x, 1)}{\partial x}  \\
\hskip-.22cm &=&\hskip-.22cm    [x^{n+1}]rx^2\frac{\partial S_r(x)}{\partial x}=rnS_n^{(r)},
\end{eqnarray*}
as desired.   \qed

\section{Further Comments}

It is possible to enumerate the aforementioned three statistics on $(a, b)$-Dyck paths by Riordan arrays, if one regards $(a, b)$-Dyck paths as colored $(a, b)$-Dyck paths such that the $\mathbf{d}$-steps in peaks colored by one of colors in $[a]$ and the other $\mathbf{d}$-steps colored by one of colors in $[b]$ for any positive integers $a$ and $b$. Let $\mathcal{C}_n(a,b)$ again be the sets of colored $(a,b)$-Dyck paths of length $2n$, the following results are derived similarly, the detailed proofs are omitted and left to the interested readers. By the relation $C(a, b; x)=S(a-b, b; x)$, Theorem \ref{theom 7.1.1}, \ref{theom 7.1.2} and \ref{theom 7.1.3} in the case $a=r$ and $b=r-1$ correspond to Theorem \ref{theom 3.1.1}, \ref{theom 4.1.1} and \ref{theom 5.1.1} respectively.

\begin{theorem}\label{theom 7.1.1}
Let $P_{n, \ell}^{(a, b)}$ denote the total number of points at level $\ell$ in all colored $(a, b)$-Dyck paths in $\mathcal{C}_{n}(a, b)$, for any integers $n\geq \ell\geq 0$, we deduce
\begin{eqnarray*}
P_{n, \ell}^{(a, b)}=\sum_{j=0}^{n-\ell}\binom{n-\ell}{j}\left\{\frac{2\ell(a-b)}{(2n-j)}\binom{2n-j}{n-\ell}+\frac{2(\ell+1)b}{2n-j+2}\binom{2n-j+2}{n-\ell} \right\}(a-b)^{j}b^{n-j-1}
\end{eqnarray*}
with $P_{0, 0}^{(a, b)}=1$. Moreover, $P_{n, \ell}^{(a, b)}+\frac{a-b}{b}\delta_{n, 0}$ is the $(n,\ell)$-entry of the Riordan array
$$\Big(C(a,b;x)^2+\frac{a-b}{b},\ bxC(a,b;x)^2\Big), $$
where $\delta_{n, 0}=1$ for $n=0$ and $0$ otherwise.
\end{theorem}

\begin{theorem}\label{theom 7.1.2}
Let $U_{n, \ell}^{(a, b)}$ denote the total number of $\mathbf{u}$-steps at level $\ell+1$ in all colored $(a, b)$-Dyck paths in $\mathcal{C}_{n+1}(a, b)$, for any integers $n\geq \ell\geq 0$, we obtain
\begin{eqnarray*}
U_{n, \ell}^{(a, b)}=\sum_{j=0}^{n-\ell}\binom{n-\ell}{j}\left\{\frac{2(\ell+1)(a-b)}{2n-j+2}\binom{2n-j+2}{n-\ell}+\frac{(2\ell+3)b}{2n-j+3}\binom{2n-j+3}{n-\ell} \right\}(a-b)^{j}b^{n-j}
\end{eqnarray*}
with $U_{0, 0}^{(a, b)}=a$. Moreover, $U_{n, \ell}^{(a, b)}$ is the $(n,\ell)$-entry of the Riordan array
$$\Big(C(a,b;x)^2\big(a-b+bC(a,b;x)\big),\ bxC(a,b;x)^2\Big). $$
\end{theorem}

\begin{theorem}\label{theom 7.1.3}
Let $p_{n, \ell}^{(a, b)}$ denote the total number of peaks at level $\ell+1$ in all colored $(a, b)$-Dyck paths in $\mathcal{C}_{n+1}(a, b)$, for any integers $n\geq \ell\geq 0$, we derive
\begin{eqnarray*}
p_{n, \ell}^{(a, b)}=\sum_{j=0}^{n-\ell}\frac{2a(\ell+1)}{2n-j+2}\binom{2n-j+2}{n-\ell}\binom{n-\ell}{j}(a-b)^{j}b^{n-j}
\end{eqnarray*}
with $p_{0, 0}^{(a, b)}=a$. Moreover, $p_{n, \ell}^{(a, b)}$ is the $(n,\ell)$-entry of the Riordan array
$$\Big(aC(a,b;x)^2,\ bxC(a,b;x)^2\Big). $$
\end{theorem}

\vskip0.5cm
\section*{Declaration of competing interest}

The authors declare that they have no known competing financial interests or personal relationships that could have
appeared to influence the work reported in this paper.

\section*{Data availability}
No data was used for the research described in the article.

\section*{Acknowledgements} {The authors are grateful to the referees for the helpful suggestions
and comments.  }

\vskip.2cm

%==============================================================================================================

\end{document}